\documentclass[12pt, a4paper, oneside]{article}
\usepackage[top=1in, bottom=1.1in, left=0.9in, right=0.5in]{geometry}
\usepackage[parfill]{parskip}
\usepackage{amsmath}
\usepackage{amsthm}
\usepackage{graphicx}
\usepackage{amssymb}
\usepackage{epstopdf}
\usepackage{setspace}
\usepackage{authblk}
\usepackage{enumerate}
\usepackage{lmodern}
\usepackage[hypcap]{caption}
\usepackage[english]{babel}
\usepackage{fancyhdr}
\addto\captionsenglish{}
\addto\captionsenglish{}
\addto\captionsenglish{}
\addto\captionsenglish{}
\addto\captionsenglish{}

\newlength\longest

\newcommand{\cM}{{\mathcal M}}

\newcommand{\R}{{\mathbb R}}

\begin{document}

\newtheorem{theorem}{Teorema}[section]
\newtheorem{lemma}[theorem]{Lemma}
\newtheorem{corollary}[theorem]{Corollary}
\newtheorem{proposition}[theorem]{Proposition}
\newtheorem{conjecture}[theorem]{Conjecture}
\newtheorem{problem}[theorem]{Problem}
\newtheorem{claim}[theorem]{Claim}
\theoremstyle{definition}
\newtheorem{assumption}[theorem]{Assumption}
\newtheorem{remark}[theorem]{Remark}
\newtheorem{definition}[theorem]{Definition}
\newtheorem{example}[theorem]{Example}
\theoremstyle{remark}
\newtheorem{notation}{Notasi}
\renewcommand{\thenotation}{}

\title{Third Version of Weak Orlicz--Morrey Spaces and Its Inclusion Properties}
\author{ Al Azhary Masta${}^1$, Siti Fatimah${}^2$ and Muhammad Taqiyuddin${}^3$}

\affil{${}^{1,2,3}$Department of Mathematics Education, Universitas Pendidikan Indonesia, Jl. Dr. Setiabudi 229, Bandung 40154\\
E-mail: ${}^{1}$alazhari.masta@upi.edu, ${}^{2}$sitifatimah@upi.edu, ${}^{3}$taqi94@hotmail.com}
\date{}

\maketitle

\begin{abstract}
Orlicz--Morrey spaces are generalizations of Orlicz spaces and Morrey spaces which were first introduced by Nakai. There are three versions of Orlicz--Morrey spaces, i.e: Nakai's (2004), Sawano--Sugano--Tanaka's (2012), and Deringoz--Guliyev--Samko's (2014) versions. On this article we will discuss the third version of weak Orlicz--Morrey space which is seen as an enlargement of third version of (strong) Orlicz--Morrey space. Similar to its first version and second version, the third version of weak Orlicz-Morrey space is considered as a generalization of weak Orlicz spaces, weak Morrey spaces, and generalized weak Morrey spaces. In this study, we will investigate some properties of the third version of weak Orlicz--Morrey spaces, especially the sufficient and necessary conditions for inclusion relations between two these spaces. One of the keys to get our result is to estimate the quasi-norm of characteristics function of open balls in $\mathbb{R}^n$.

\bigskip

\noindent{\bf Keywords}: Weak Orlicz spaces, Weak Morrey spaces, Weak Orlicz-Morrey space of third version, Inclusion property.\\
{\textbf{AMS Subject Classification}}: Primary 46E30; Secondary 46B25, 42B35.
\end{abstract}

\section{Introduction}

Orlicz-Morrey spaces are generalization of Orlicz spaces and Morrey spaces and it is firstly introduced by E. Nakai in 2004 \cite{Gala,Nakai1,Nakai2}.  These spaces are one of the important topics in mathematical analysis, particularly in harmonic analysis. There are three versions of Orlicz--Morrey spaces, i.e: Nakai's (2004), Sawano--Sugano--Tanaka's (2012) \cite{Gala}, and Deringoz--Guliyev--Samko's (2014) \cite{Deringoz,Guliyev} versions. 


For a Young function $\Theta: [ 0, \infty ) \rightarrow [ 0, \infty ) $ (i.e.  $\Theta$ is convex, $\lim \limits_{t\to 0}\Theta(t)=0=\Theta(0)$, 
continuous and $ \lim \limits_{t\to\infty} \Theta(t) = \infty$), we define $\Theta^{-1}(s) :=\inf \{ r \geq 0 : \Theta (r) > s \}$. Given two Young functions $\Theta_1, \Theta_2$, we write $\Theta_1 \prec \Theta_2$ if there exists a constant $C > 0$ such that $\Theta_1(t) \leq \Theta_2(Ct)$ for all $ t > 0$.

Now, let $G_{\theta}$ be the set of all functions $\theta : (0, \infty) \rightarrow (0, \infty)$
such that $\theta(r)$ is decreasing but $\Theta^{-1}(t^{-n})\theta(t)^{-1}$ is almost decreasing for all $ t>0$. Let $\theta_1 \in G_{\theta_1} $ and $\theta_2 \in G_{\theta_2}$, we denote $\theta_1 \lesssim \theta_2$ if there exists a constant $C > 0$ such that $\theta_1(t)  \leq C \theta_2(t)$ for all $ t >0$.

First we recall definition of (strong) Orlicz--Morrey spaces of Deringoz--Guliyev--Samko's (2014) version. Let $\Theta$ be a Young function and $\theta \in G_{\theta}$, 
\textit{the Orlicz--Morrey space} $\cM_{\theta,\Theta}(\mathbb{R}^{n})$ is the set of measurable
functions $f$ on $\R^n$ such that
$$\| f \|_{\cM_{\theta,\Theta}(\mathbb{R}^{n})} := \sup\limits_{a \in \mathbb{R}^n,~r>0} \frac{1}{\theta(|B(a,r)|^{\frac{1}{n}})}\Theta^{-1}\Bigl(\frac{1}{|B(a,r)|}\Bigr) \|f\|_{L_{\Theta}(B(a,r))} < \infty, $$
where $\| f \|_{L_{\Theta}(B(a,r))} := \inf \bigl \{  {b>0:
	\int_{B(a,r)}\Theta \bigl(\frac{|f(x)|}{b} \bigr) dx \leq1}\bigr \}$.
Here, $B(a,r)$ denotes the open ball in $\mathbb{R}^n$ centered at $a \in \R^n$ with radius $r > 0$, and $|B(a,r)|$ for its Lebesgue measure.

Meanwhile, for $\Theta$ is Young function and $\theta \in G_{\theta}$, \textit{the weak Orlicz--Morrey space} $w\cM_{\theta,\Theta}(\mathbb{R}^{n})$ is the set of all measurable functions $f$ on $\R^n$ such that 
 
 $$\| f \|_{w\cM_{\theta,\Theta}(\mathbb{R}^{n})} := \sup\limits_{a \in
 	\mathbb{R}^n,~r>0} \frac{1}{\theta(|B(a,r)|^{\frac{1}{n}})}\Theta^{-1}\Bigl(\frac{1}{|B(a,r)|}\Bigr) \|f\|_{wL_{\Theta}(B(a,r))} < \infty, $$
 where $\| f \|_{wL_{\Theta}(B(a,r))} := \inf \left\{  {b>0:
 	\mathop {\sup }\limits_{t > 0} \Phi(t) \Bigl| \Bigl\{ x \in B(a,r) : \frac{|f(x)|}{b} > t \Bigr\}
 	\Bigr| \leq1}\right\}.$

The space $w\cM_{\theta,\Theta}(\mathbb{R}^{n})$ is quasi-Banach spaces equipped with the quasi-norm $\| \cdot \|_{w\cM_{\theta,\Theta}(\mathbb{R}^{n})}$. Note that, analog with $\cM_{\theta,\Theta}(\mathbb{R}^{n})$ space, $w\cM_{\theta,\Theta}(\mathbb{R}^{n})$ also covers many classical spaces, which shown in the following example.

\bigskip

\begin{example}
Let $1 \leq p \leq q < \infty$, $\Phi$ be a Young function, and $\theta \in G_{\theta}$ then we obtain:

\begin{enumerate}
\item If $\Theta(t)=t^p$ and $\theta(t)=t^{\frac{-n}{p}}$, then $w\cM_{\theta,\Theta}(\mathbb{R}^{n}) = wL^p(\mathbb{R}^{n})$ is weak Lebesgue space. 

\item If $\Theta(t)=t^q$ and $\theta(t)=t^{\frac{-n}{p}}$, then $w\cM_{\theta,\Theta}(\mathbb{R}^{n}) = w\cM_{p}^q(\mathbb{R}^{n})$ is classical weak Morrey space.

\item If $\Theta(t)=t^p$, then $w\cM_{\theta,\Theta}(\mathbb{R}^{n}) = w\cM_{\theta}^p (\mathbb{R}^{n})$ is generalized weak Morrey space.

\item  If $\theta(t)=\Theta^{-1}(t^{-n})$, then $w\cM_{\theta,\Theta}(\mathbb{R}^{n}) = wL_{\Theta}(\mathbb{R}^{n})$ is weak Orlicz space. 
\end{enumerate}	
\end{example}

Moreover, the relationship between $\cM_{\theta,\Theta}(\mathbb{R}^{n})$ space and $w\cM_{\theta,\Theta}(\mathbb{R}^{n})$ space can be stated as the following lemma.

\bigskip
\begin{lemma}\label{lemma:1.2}
	Let $\Theta$ be a Young function and $\theta \in G_{\theta}$. Then $\cM_{\theta,\Theta}(\mathbb{R}^{n})
	\subseteq w\cM_{\theta,\Theta}(\mathbb{R}^{n})$ with $\| f \|_{w\cM_{\theta,\Theta}(\mathbb{R}^{n})} \leq \| f \|_{\cM_{\theta,\Theta}(\mathbb{R}^{n})}$ for every $f \in \cM_{\theta,\Theta}(\mathbb{R}^{n})$.
\end{lemma}

\bigskip
 Many authors have been culminating important observations about inclusion properties of function spaces, see \cite{Kufner,Lech,Masta1,Masta2,Masta3,Taqi}, etc. Recently, Masta \textit{et al.} \cite{Masta3} obtained sufficient and necessary conditions for inclusion of (strong) Orlicz--Morrey spaces of all versions. In the same paper, Masta \textit{et al.} also proved the sufficient and necessary conditions for inclusion properties of weak Orlicz--Morrey spaces of Nakai's and Sawano--Sugano--Tanaka's versions.
  
 In this paper, we would like to obtain the inclusion properties of weak Orlicz--Morrey space $w\cM_{\theta,\Theta}(\mathbb{R}^{n})$ of Deringoz--Guliyev--Samko's version, and compare it with the result for Nakai's and Sawano-Sugano-Tanaka's versions.

To achieve our purpose, we will use the similar methods in \cite{Gunawan, Masta1, Masta2, Alen} which pay attention to the characteristic functions of open balls in $\R^n$, in the following lemma.

\bigskip
\begin{lemma}\label{lemma:1.3} \cite{Guliyev}
Let $\Theta$ be a Young function, $\theta \in G_{\theta}$, and $ r_0>0$, then there exists a constant $C>0$ such that $$ \frac{1}{\theta(r_0)} \leq \| \chi_{B(0,r_0)} \|_{\cM_{\theta,\Theta}(\mathbb{R}^{n})} \leq \frac{C}{\theta(r_0)}$$.
\end{lemma}

\bigskip
For weak Orlicz--Morrey spaces, we have the following lemma.

\bigskip
\begin{lemma}\label{lemma:1.4}
	Let $\Theta$ be a Young function, $\theta \in G_{\theta}$, and $ r_0>0$, then there exists a constant $C>0$ such that $$ \frac{1}{\theta(r_0)} \leq \| \chi_{B(0,r_0)} \|_{w\cM_{\theta,\Theta}(\mathbb{R}^{n})} \leq \frac{C}{\theta(r_0)}$$.
\end{lemma}

Proof. Since $\Theta$ is a Young function and $\theta \in G_{\theta}$, by Lemmas \ref{lemma:1.2} and \ref{lemma:1.3},
we have $$\| \chi_{B(0,r_{0})} \|_{w\cM_{\theta,\Theta}(\mathbb{R}^{n})} \leq \| \chi_{B(0,r_{0})} \|_{\cM_{\theta,\Theta}(\mathbb{R}^{n})} \leq
\frac{C}{\theta(r_0)}.$$ On the other hand,
\begin{align*}
\| \chi_{B(0,r_{0})} \|_{w\cM_{\theta,\Theta}(\mathbb{R}^{n})}
& = \sup\limits_{a \in \mathbb{R}^n, r > 0} \frac{1}{\theta(|B(a,r)|^{\frac{1}{n}})}\Theta^{-1}\Bigl(\frac{1}{|B(a,r)|}\Bigr) \| \chi_{B(0,r_{0})} \|_{wL_{B(a,r)}}\\
& = \sup\limits_{a \in \mathbb{R}^n, r > 0} \frac{1}{\theta(|B(a,r)|^{\frac{1}{n}})}\Theta^{-1}\Bigl(\frac{1}{|B(a,r)|}\Bigr) \frac{1}{\Theta^{-1} \bigl( \frac{|B(a,r)|}{|B(a,r)
		\cap B(0,r_{0})|} \bigr)}\\
& \geq \frac{1}{\theta (r_0)}.
\end{align*}

Consequently, we have $\frac{1}{\theta(r_0)} \leq \| \chi_{B(0,r_0)} \|_{w\cM_{\theta,\Theta}(\mathbb{R}^{n})} \leq \frac{C}{\theta(r_0)}$. \qed

\medskip

In this paper, the letter $C$ will be used for constants that may change from line to line, while constants with subscripts, such as $C_{1}, C_{2}$, do not change in different lines.

\section{Results}

First, we reprove sufficient and necessary conditions for inclusion properties of Orlicz--Morrey space $\cM_{\theta, \Theta}(\mathbb{R}^n)$ in the following theorem.

\medskip

\begin{theorem}\label{theorem:2.1} \cite{Masta3}
	Let $\Theta_1, \Theta_2$ be Young functions such that $\Theta_1 \prec \Theta_2$, $\Theta^{-1}_1 \lesssim \Theta^{-1}_2$, $\theta_1\in G_{\theta_1}$ and $\theta_2 \in G_{\theta_2}$. Then the following statements are equivalent:

	{\parindent=0cm
		{\rm (1)} $\theta_2 \lesssim \theta_1$.
		
		{\rm (2)} $\cM_{\theta_2, \Theta_2}(\mathbb{R}^n) \subseteq \cM_{\theta_1, \Theta_1}(\mathbb{R}^n)$.
		
		{\rm(3)} There exists a constant $C > 0$ such that
		$$\| f \|_{\cM_{\theta_1, \Theta_1}(\mathbb{R}^n)} \leq C \| f \|_{\cM_{\theta_2, \Theta_2} (\mathbb{R}^n)}$$
		for every $ f \in \cM_{\theta_2, \Theta_2}(\mathbb{R}^n)$.
		\par}
\end{theorem}

\textit{Proof.} Assume that (1) holds and let $ f \in \cM_{\theta_2, \Theta_2}(\mathbb{R}^n)$.
Since $\Theta_1 \prec \Theta_2$, by using similar arguments in the proof of Corollary 2.3 in \cite{Masta1}, we have

$$\|f\|_{L_{\Theta_1}(B(a,r))} \leq C \|f\|_{L_{\Theta_2}(B(a,r))}$$

for every $ B(a,r) \subseteq \R^n$.

Knowing that, $\Theta^{-1}_1 \lesssim \Theta^{-1}_2$ and $\theta_1 \lesssim \theta_2$ ( i.e. there exists constant
$C_1,C_2 > 0$ such that $\Theta^{-1}_1(t) \leq C_1\Theta^{-1}_2(t)$ and
$\theta_2(t)  \leq C_2 \theta_1(t)$ for every $ t > 0 $), we obtain
\begin{align*}
\| f \|_{\cM_{\theta_1,\Theta_1}(\mathbb{R}^{n})} &:= \sup\limits_{a \in
	\mathbb{R}^n,~r>0} \frac{1}{\theta_1(|B(a,r)|^{\frac{1}{n}})}\Theta_1^{-1}\Bigl(\frac{1}{|B(a,r)|}\Bigr) \|f\|_{L_{\Theta_1}(B(a,r))}\\
&\leq \sup\limits_{a \in
	\mathbb{R}^n,~r>0} \frac{C}{\theta_1(|B(a,r)|^{\frac{1}{n}})}\Theta_1^{-1}\Bigl(\frac{1}{|B(a,r)|}\Bigr) \|f\|_{L_{\Theta_2}(B(a,r))}\\
&\leq \sup\limits_{a \in
	\mathbb{R}^n,~r>0} \frac{CC_1}{\theta_1(|B(a,r)|^{\frac{1}{n}})}\Theta_2^{-1}\Bigl(\frac{1}{|B(a,r)|}\Bigr) \|f\|_{L_{\Theta_2}(B(a,r))}\\
&\leq \sup\limits_{a \in
	\mathbb{R}^n,~r>0} \frac{CC_1C_2}{\theta_2(|B(a,r)|^{\frac{1}{n}})}\Theta_2^{-1}\Bigl(\frac{1}{|B(a,r)|}\Bigr) \|f\|_{L_{\Theta_2}(B(a,r))}\\
&:=CC_1C_2\| f \|_{\cM_{\theta_2,\Theta_2}(\mathbb{R}^{n})}
\end{align*}

This proves that $\cM_{\theta_2, \Theta_2}(\mathbb{R}^{n}) \subseteq \cM_{\theta_1, \Theta_1} (\mathbb{R}^{n})$.

Next, since $(\cM_{\theta_2, \Theta_2}(\mathbb{R}^n), \cM_{\theta_1, \Theta_1}(\mathbb{R}^n))$ is a Banach pair,
it follows from \cite[Lemma 3.3]{Krein} that (2) and (3) are equivalent. It thus remains to show that (3) implies (1).

Assume that (3) holds. Let $ r_0 > 0$. By Lemma \ref{lemma:1.3}, we have
$$
\frac{1}{\theta_1(r_0)} \leq \| \chi_{B(0,r_{0})} \|_{\cM_{\theta_1, \Theta_1}(\mathbb{R}^n)} \leq C \| \chi_{B(0,r_{0})} \|_{\cM_{\theta_2, \Theta_2}(\mathbb{R}^n)} \leq \frac{C}{\theta_2(r_0)}, $$

Since $r_{0} > 0$ is arbitrary, we conclude that $\theta_2(t) \leq C\theta_1(t)$ for every $t > 0$. \qed

Now we come to the inclusion property of weak Orlicz--Morrey spaces $w\cM_{\theta_1, \Theta_1}(\mathbb{R}^n)$ and $w\cM_{\theta_2, \Theta_2}(\mathbb{R}^n)$ with respect to Young functions $\Theta_1, \Theta_2$ and parameters $\theta_1, \theta_2$.
\medskip
  
\begin{theorem}\label{theorem:2.2}
	Let $\Theta_1, \Theta_2$ be Young functions such that $\Theta_1 \prec \Theta_2$, $\Theta^{-1}_1 \lesssim \Theta^{-1}_2$, $\theta_1\in G_{\theta_1}$ and $\theta_2 \in G_{\theta_2}$. Then the following statements are equivalent:

	{\parindent=0cm
		{\rm (1)} $\theta_2 \lesssim \theta_1$.
		
		{\rm (2)} $w\cM_{\theta_2, \Theta_2}(\mathbb{R}^n) \subseteq w\cM_{\theta_1, \Theta_1}(\mathbb{R}^n)$.
		
		{\rm(3)} There exists a constant $C > 0$ such that
		$$\| f \|_{w\cM_{\theta_1, \Theta_1}(\mathbb{R}^n)} \leq C \| f \|_{w\cM_{\theta_2, \Theta_2} (\mathbb{R}^n)}$$
		for every $ f \in w\cM_{\theta_2, \Theta_2}(\mathbb{R}^n)$.
		\par}
\end{theorem}

\textit{Proof.} Assume that (1) holds and let $ f \in w\cM_{\theta_2, \Theta_2}(\mathbb{R}^n)$.
Since $\Theta_1 \prec \Theta_2$, by using similar arguments in the proof of Theorem 3.3 in \cite{Masta1}, we have

$$\|f\|_{wL_{\Theta_1}(B(a,r))} \leq C \|f\|_{wL_{\Theta_2}(B(a,r))}$$

for every $ B(a,r) \subseteq \R^n$.

Knowing that, $\Theta^{-1}_1 \lesssim \Theta^{-1}_2$ and $\theta_1 \lesssim \theta_2$ (i.e there exists constant
$C_1,C_2> 0$ such that $\Theta^{-1}_1(t) \leq C_1\Theta^{-1}_2(t)$ and
$\theta_2(t)  \leq C_2 \theta_1(t)$ for every $ t > 0 $), we obtain
\begin{align*}
\| f \|_{w\cM_{\theta_1,\Theta_1}(\mathbb{R}^{n})} &:= \sup\limits_{a \in
	\mathbb{R}^n,~r>0} \frac{1}{\theta_1(|B(a,r)|^{\frac{1}{n}})}\Theta_1^{-1}\Bigl(\frac{1}{|B(a,r)|}\Bigr) \|f\|_{wL_{\Theta_1}(B(a,r))}\\
&\leq \sup\limits_{a \in
	\mathbb{R}^n,~r>0} \frac{C}{\theta_1(|B(a,r)|^{\frac{1}{n}})}\Theta_1^{-1}\Bigl(\frac{1}{|B(a,r)|}\Bigr) \|f\|_{wL_{\Theta_2}(B(a,r))}\\
&\leq \sup\limits_{a \in
	\mathbb{R}^n,~r>0} \frac{CC_1}{\theta_1(|B(a,r)|^{\frac{1}{n}})}\Theta_2^{-1}\Bigl(\frac{1}{|B(a,r)|}\Bigr) \|f\|_{wL_{\Theta_2}(B(a,r))}\\
&\leq \sup\limits_{a \in
	\mathbb{R}^n,~r>0} \frac{CC_1C_2}{\theta_2(|B(a,r)|^{\frac{1}{n}})}\Theta_2^{-1}\Bigl(\frac{1}{|B(a,r)|}\Bigr) \|f\|_{wL_{\Theta_2}(B(a,r))}\\
&:=CC_1C_2\| f \|_{w\cM_{\theta_2,\Theta_2}(\mathbb{R}^{n})}
\end{align*}

This proves that $w\cM_{\theta_2, \Theta_2}(\mathbb{R}^{n}) \subseteq w\cM_{\theta_1, \Theta_1} (\mathbb{R}^{n})$.

As a corollary of the Open Mapping Theorem [3, Appendix G], we are aware that [6, Chapter I,
Lemma 3.3] still holds for quasi-Banach spaces, and so (2) and (3) are equivalent.

Assume that (3) holds. Let $ r_0 > 0$. By Lemma \ref{lemma:1.4}, we have
$$
\frac{1}{\theta_1(r_0)} \leq \|\chi_{B(0,r_{0})} \|_{w\cM_{\theta_1, \Theta_1}(\mathbb{R}^n)} \leq C \| \chi_{B(0,r_{0})} \|_{w\cM_{\theta_2, \Theta_2}(\mathbb{R}^n)} \leq \frac{C}{\theta_2(r_0)}, $$

Since $r_{0} > 0$ is arbitrary, we conclude that $\theta_2(t) \leq C\theta_1(t)$ for every $t > 0$. \qed

For generalized weak Morrey spaces, we have the following corollary.

\begin{corollary}
	Let $1 \leq p < \infty$, $\theta_1\in G_{\theta_1}$ and $\theta_2 \in G_{\theta_2}$. Then the following statements are equivalent:

{\parindent=0cm
	{\rm (1)} $\theta_2 \lesssim \theta_1$.
	
	{\rm (2)} $w\cM_{\theta_2}^p(\mathbb{R}^n) \subseteq w\cM_{\theta_1}^p(\mathbb{R}^n)$.
	
	{\rm(3)} There exists a constant $C > 0$ such that
	$$\| f \|_{w\cM_{\theta_1}^p(\mathbb{R}^n)} \leq C \| f \|_{w\cM_{\theta_2}^p (\mathbb{R}^n)}$$
	for every $ f \in w\cM_{\theta_2}^p(\mathbb{R}^n)$.
	\par}
\end{corollary}
\section{Concluding Remarks}

We have shown the sufficient and necessary conditions for the inclusion relation between weak Orlicz--Morrey space $w\cM_{\theta, \Theta} (\R^n)$. In the proof of the inclusion property we used the norm of characteristic function on $\R^n$. The inclusion properties of weak Orlicz-Morrey space $w\cM_{\theta, \Theta} (\R^n)$ (Theorem \ref{theorem:2.2}) and weak Orlicz--Morrey space $w\cM_{\psi, \Psi} (\R^n)$ of Sawano--Sugano--Tanaka's version (Theorem 3.9 in \cite{Masta3}) generalize the inclusion properties of weak Morrey spaces and generalized weak Morrey spaces in \cite{Gunawan}. Meanwhile, the inclusion properties of weak Orlicz-Morrey space $wL_{\phi,\Phi} (\R^n)$ of Nakai's version (Theorem 3.4 in \cite{Masta3}) generalize the inclusion properties of weak Orlicz spaces in \cite{Masta1}.

Comparing Theorem 2.2 and Theorem 3.9 in \cite{Masta3}, we can say that the condition on the Young function for the inclusion of the weak Orlicz--Morrey space $w\cM_{\psi, \Psi} (\R^n)$ is simpler than that for the weak Orlicz--Morrey space $w\cM_{\theta, \Theta} (\R^n)$.

As a corollary of Lemma \ref{lemma:1.2}, Theorem \ref{theorem:2.1} and Theorem \ref{theorem:2.2}, we also have the following inclusion
relations
\[
\begin{array}{ccc}
\cM_{\theta_2, \Theta_2} & \rightarrow &\cM_{\theta_1, \Theta_1}\\
\downarrow & \searrow & \downarrow\\
w\cM_{\theta_2, \Theta_2} &\rightarrow & w\cM_{\theta_1, \Theta_1}
\end{array}
\]
for $\Theta_1 \prec \Theta_2$, $\Theta^{-1}_1 \lesssim \Theta^{-1}_2$ and $\theta_2 \lesssim \theta_1$, where the arrows
mean `contained in' or `embedded into'.


\begin{thebibliography}{10}

\bibitem{Deringoz}
F. Deringoz, V. S. Guliyev, and S. Samko, ``Boundedness of the maximal and singular operators on generalized Orlicz--Morrey spaces'', In: Operator Theory, Operator Algebra and Applications,
\emph{Operator Theory: Advances and Applications} \textbf{242} (2014), 139--158.

\bibitem{Gala}
S. Gala, Y. Sawano, and H. Tanaka, ``A remark on two generalized Orlicz-Morrey spaces'', \emph{J. Approx. Theory} \textbf{198} (2015), 1--9.

\bibitem{Gra09}
L.~Grafakos, \textit{Modern Fourier analysis}, Graduate Texts in Mathematics, {\bf 250}, Springer, New York, 2009.

\bibitem{Guliyev}
V. S. Guliyev, S. G. Hasanov, Y. Sawano, and T. Noi, ``Non-smooth atomic decompositions for generalized Orlicz--Morrey spaces of the third kind'', \emph{Acta Appl. Math.} {\bf 145}--1 (2017), 133--174 [DOI: 10.1007/s10440-016-0052-7].


\bibitem {Gunawan}
H. Gunawan, D.I. Hakim, K.M. Limanta, A.A. Masta, ``Inclusion properties of generalized Morrey spaces'', \emph{Math. Nachr.} \textbf{290} (2017), 332-340 [DOI:10.1002/mana.201500425].


\bibitem {Krein}
S.G Kre\v{i}n, Yu.\={I} Petun\={i}n, and E.M. Sem\"{e}nov, \emph{Interpolation of Linear
Operators}, Translation of Mathematical Monograph vol. 54, American Mathematical Society, Providence, R.I., 1982.

\bibitem{Kufner}
A. Kufner, O. John, and S. Fu\"{c}ik, \emph{Function Spaces}, Noordhoff International
Publishing, Czechoslovakia, 1977.


\bibitem{Lech} L. Maligranda, \emph{Orlicz Spaces and Interpolation}, Departamento de Matem\'{a}tica, Universidade Estadual de Campinas, 1989.

\bibitem{Masta1} A.A. Masta, H. Gunawan, and W. Setya-Budhi, ``Inclusion property of Orlicz and weak Orlicz spaces'', \emph{J. Math. Fund. Sci}. \textbf{48}-3 (2016), 193--203 [DOI: http://dx.doi.org/10.56142Fj.math.fund.sci.2016.48.3.1].

\bibitem{Masta2} A.A. Masta, H. Gunawan, W.S. Budhi, ``An inclusion property of Orlicz-Morrey spaces'', \emph{J. Phys.: Conf. Ser.}\textbf{893 }012015 (2017), 1--7, [DOI: https://doi.org/10.1088/1742-6596/893/1/012015].

\bibitem{Masta3} A.A. Masta, H. Gunawan, and W. Setya-Budhi, ``On inclusion properties of two versions of Orlicz--Morrey
spaces'', \emph{Mediterr. J. Math.}, \textbf{14}-6 (2017), 1--12 [DOI: https://doi.org/10.1007/s00009-017-1030-7].

\bibitem{Nakai1}	
E. Nakai, ``On Orlicz-Morrey spaces'', research report
[http://repository.kulib.kyoto-u.ac.jp/dspace/bitstream/2433/58769/1/1520-
10.pdf,
accessed on August 17, 2015.]

\bibitem{Nakai2}
E. Nakai, ``Orlicz-Morrey spaces and some integral operators'', research report
[http://repository.kulib.kyoto-u.ac.jp/dspace/bitstream/2433/26035/1/
1399-13.pdf,
accessed on August 17, 2015.]


\bibitem{Alen}
A. Osan\c{c}liol, ``Inclusion between weighted Orlicz spaces'', \emph{J. Inequal. Appl.}
\textbf{2014}-390 (2014), 1--8, [DOI: https://doi.org/10.1186/1029-242X-2014-390].

\bibitem{Taqi} M. Taqiyuddin M and A. A.Masta, ``Inclusion properties of Orlicz spaces and weak Orlicz spaces generated by concave function'', \textit{IOP Conf. Ser.: Mater. Sci. Eng.,} \textbf{288}-012103, 1--5 [DOI: https://doi.org/10.1088/1757-899X/288/1/012103].


\end{thebibliography}
\end{document}